\documentclass[11pt,reqno]{amsart}
\usepackage{amsmath,amsthm,amssymb,verbatim,enumerate,ifthen,times}
\usepackage[mathscr]{eucal}
\usepackage[utf8]{inputenc}
\usepackage[T1]{fontenc}
\newcommand{\N}{\mathbb{N}}

\newcommand{\R}{\mathbb{R}}
\newcommand{\C}{\mathscr{C}}
\newcommand{\CC}{\mathbb{C}}
\newcommand{\A}{\mathbb{A}}

\newcommand{\B}{\mathbb{B}}

\DeclareMathOperator{\sgn}{sign}

\newtheorem{theorem}{Theorem}
\newtheorem*{theorem*}{Theorem}
\def\Thm#1#2{\ifthenelse{\equal{#1}{*}}{\begin{theorem*}#2\end{theorem*}}
  {\begin{theorem}\label{T#1}#2\end{theorem}}}
\newtheorem{Atheorem}{Theorem}

\def\thm#1{Theorem~\ref{T#1}}

\newtheorem{proposition}[theorem]{Proposition}
\newtheorem*{proposition*}{Proposition}
\def\Prp#1#2{\ifthenelse{\equal{#1}{*}}{\begin{proposition*}#2\end{proposition*}}
             {\begin{proposition}\label{P#1}#2\end{proposition}}}

\newtheorem{corollary}[theorem]{Corollary}
\newtheorem*{corollary*}{Corollary}
\def\Cor#1#2{\ifthenelse{\equal{#1}{*}}{\begin{corollary*}#2\end{corollary*}}
             {\begin{corollary}\label{C#1}#2\end{corollary}}}

\newtheorem{lemma}[theorem]{Lemma}
\newtheorem*{lemma*}{Lemma}
\def\Lem#1#2{\ifthenelse{\equal{#1}{*}}{\begin{lemma*}#2\end{lemma*}}
             {\begin{lemma}\label{L#1}#2\end{lemma}}}
\def\lem#1{Lemma~\ref{L#1}}

\newtheorem{example}[theorem]{Example}
\newtheorem*{example*}{Example}
\def\Exa#1#2{\ifthenelse{\equal{#1}{*}}{\begin{example*}\rm #2\end{example*}}
             {\begin{example}\label{Ex#1}\rm #2\end{example}}}

\newtheorem{problem}[theorem]{Problem}

\theoremstyle{definition}
\newtheorem{definition}[theorem]{Definition}

\newtheorem{remark}[theorem]{Remark}
\newtheorem*{remark*}{Remark}
\def\Rem#1#2{\ifthenelse{\equal{#1}{*}}{\begin{remark*}\rm #2\end{remark*}}
             {\begin{remark}\label{R#1}\rm #2\end{remark}}}

\newcommand{\eq}[1]{\eqref{E#1}}
\newcommand{\Eq}[2]{\ifthenelse{\equal{#1}{*}}
  {\begin{equation*}\begin{aligned}[]#2\end{aligned}\end{equation*}}
  {\begin{equation}\begin{aligned}[]\label{E#1}#2\end{aligned}\end{equation}}}

\long\def\comment#1{}

\begin{document}
\large

\date{\today}

\title[Equality of Cauchy means and quasiarithmetic means]
{Characterization of the equality of Cauchy means \\ to quasiarithmetic means}

\author[R. L. Lovas]{Rezső L. Lovas}
\author[Zs. P\'ales]{Zsolt P\'ales}
\address[R. L. Lovas and Zs. P\'ales]{Institute of Mathematics, University of Debrecen, H-4002 Debrecen, Pf.\ 400, Hungary}
\email{\{lovas,pales\}@science.unideb.hu}

\author[A. Zakaria]{Amr Zakaria}
\address[A. Zakaria]{Doctoral School of Mathematical and Computational Sciences, University of Debrecen, H-4002 Debrecen, Pf.\ 400, Hungary; Department of Mathematics, Faculty of Education, Ain Shams University, Cairo 11341, Egypt}
\email{amr.zakaria@edu.asu.edu.eg}

\subjclass[2010]{39B40, 26E60}
\keywords{Cauchy mean; quasiarithmetic mean; equality problem; functional equation}


\thanks{The research of the second author was supported by the EFOP-3.6.1-16-2016-00022 project. This project is co-financed by the European Union and the European Social Fund.}

\begin{abstract}
The main result of this paper provides six necessary and sufficient conditions under various regularity assumptions for a so-called Cauchy mean to be identical to a two-variable quasiarithmetic mean. One of these conditions says that a Cauchy mean is quasiarithmetic if and only if the range of its generating functions is covered by a nondegenerate conic section.
\end{abstract}
\maketitle

\section{Introduction}

Throughout this paper, the symbols $\R$, $\R_+$, and $\N$ will denote the sets of real, positive real, and natural numbers, respectively, and $I$ will always denote a nonempty open interval. Given a continuous strictly monotone function $\varphi:I\to\R$, the \emph{two-variable quasiarithmetic mean} $\A_\varphi:I^2\to I$ is defined by
\Eq{*}{
\A_\varphi(x,y):=\varphi^{-1}\left(\frac{\varphi(x)+\varphi(y)}{2}\right).
}
A systematic study of these means can be found in the book \cite{HarLitPol34}. 
A characterization theorem of these means was obtained by Aczél in \cite{Acz48a} (cf. also \cite{Acz66}, \cite{AczDho89}).

There are two essential generalizations of two-variable quasiarithmetic means. The first one is due to Bajraktarevi\'c \cite{Baj58}: Given two functions $f,g:I\to\R$ such that $g$ is nowhere zero and $f/g$ is strictly monotone and continuous, the \emph{two-variable Bajraktarevi\'c mean} $\B_{f,g}:I^2\to I$ is defined by
\Eq{*}{
  \B_{f,g}(x,y)
  :=\bigg(\frac{f}{g}\bigg)^{-1}\left(\frac{f(x)+f(y)}{g(x)+g(y)}\right).
}
It is immediate to see that $\B_{\varphi,1}\equiv \A_\varphi$, showing that two-variable quasiarithmetic means form a proper subclass of two-variable Bajraktarevi\'c means.

The second generalization is due to Leach and Sholander \cite{LeaSho84} (cf.\ also Losonczi \cite{Los00a}):
Given two differentiable functions $f,g:I\to\R$ such that $g'$ is nowhere zero, $f'/g'$ is strictly monotone, the \emph{Cauchy mean} $\CC_{f,g}:I^2\to I$ is defined by  
\Eq{*}{
\CC_{f,g}(x,y)
:=\begin{cases}
  \bigg(\dfrac{f'}{g'}\bigg)^{\!\!-1}\!\!\bigg(\dfrac{f(x)-f(y)}{g(x)-g(y)}\bigg) 
  &\mbox{if } x\neq y,\\
  x&\mbox{if } x=y.
  \end{cases}
}
Observe that if $\varphi$ is differentiable with a nonvanishing derivative, then
$\CC_{\varphi^2,\varphi}\equiv \A_\varphi$, i.e., Cauchy means contain the class of two-variable quasiarithmetic means with a differentiable generating function whose derivative is nonvanishing.

In the sequel, we say that two pairs of functions $(f,g):I\to\R^2$ and $(h,k):I\to\R^2$ are equivalent (and we write $(f,g)\sim(h,k)$) if there exist real constants $a,b,c,d$ with $ad\neq bc$ such that
\Eq{abcd}{
  h=af+bg \qquad\mbox{and}\qquad k=cf+dg.
}
One can easily check that $\sim$ is an equivalence relation, indeed. 

For a real parameter $p\in\R$, we introduce the sine and cosine type functions $S_p,C_p:\R\to\R$ by
\Eq{*}{
  S_p(x):=\begin{cases}
           \sin(\sqrt{-p}x) & \mbox{ if } p<0, \\
           x & \mbox{ if } p=0, \\
           \sinh(\sqrt{p}x) & \mbox{ if } p>0, 
         \end{cases}\qquad\mbox{and}\qquad
  C_p(x):=\begin{cases}
           \cos(\sqrt{-p}x) & \mbox{ if } p<0, \\
           1 & \mbox{ if } p=0, \\
           \cosh(\sqrt{p}x) & \mbox{ if } p>0. 
        \end{cases}
}
It is easily seen that the functions $S_p$ and $C_p$ form a fundamental system of solutions for the second-order homogeneous linear differential equation $Y''=pY$.

We introduce the following regularity classes for the generating functions of Cauchy means as follows: Let the class $\C_1(I)$ contain all pairs $(f,g)$ such that
\begin{enumerate}
 \item[(i)] $f,g:I\to\R$ are continuously differentiable functions such that $g'$ is nowhere zero on $I$.
 \item[(ii)] $f'/g'$ is strictly monotone on $I$.
\end{enumerate}
For $n\geq2$, let $\C_n(I)$ denote the class of all pairs $(f,g)$ such that
\begin{enumerate}
 \item[(+i)] $f,g:I\to\R$ are $n$ times continuously differentiable functions and $g'$ is nowhere zero on $I$.
 \item[(+ii)] $(f'/g')'$ is nowhere zero on $I$.
\end{enumerate}

Finally, for $(f,g)\in\C_n(I)$ and for $i,j\in\{0,\dots,n\}$, we define the generalized Wronski-type determinant $W^{i,j}_{f,g}:I\to\R$ by
\Eq{W}{
  W^{i,j}_{f,g}:=\left|\begin{array}{cc} f^{(i)} & f^{(j)} \\ g^{(i)} & g^{(j)} \end{array}\right|.
}

The equality of Cauchy means to two-variable quasiarithmetic means has been characterized by Kiss and Páles \cite{KisPal19}.

\Thm{0}{
Let $(f,g)\in\C_1(I)$ and $h:I\to\R$ be continuous and strictly monotone. Then 
\Eq{*}{
\CC_{f,g}(x,y)=\A_h(x,y) \qquad(x,y\in I)
}
holds if and only if $h$ is differentiable with a nonvanishing first derivative and there exists a constant $p\in\R$ such that 
\Eq{equ}{
(f',g')\sim(h'\cdot S_p\circ h,h'\cdot C_p\circ h).
}}

The following theorem provides characterization of the equality of two-variable Bajraktarevi\'c means to two-variable quasiarithmetic means (cf.\ \cite{DarMakPal04}, \cite{KisPal18}, \cite{PalZak19c}).

\Thm{0+}{
Let $f,g:I\to\R$ be two functions such that $g$ is everywhere positive on $I$ and $f/g$ is strictly monotone and continuous on $I$. Then the following statements are equivalent.
\begin{enumerate}[(i)]
\item There exists a continuous strictly monotone function $h:I\to\R$ such that
\Eq{fgh}{
  \B_{f,g}(x,y)=\A_h(x,y)\qquad(x,y\in I).
}
\item There exist a continuous strictly monotone function $h:I\to\R$ and a constant $p\in\R$ such that
\Eq{gf}{
  (f,g)\sim(S_p\circ h,C_p\circ h).
}
\item There exist real constants $\alpha,\beta,\gamma$ such that
\Eq{abc}{
  \alpha f^2+\beta fg+\gamma g^2=1.
}
\item Provided that $f$ and $g$ are continuously differentiable and $W^{1,0}_{f,g}$ is nowhere zero on $I$, equation \eq{fgh}
holds with $h=\int W^{1,0}_{f,g}$.
\item Provided that $f$ and $g$ are twice continuously differentiable and $W^{1,0}_{f,g}$ is nowhere zero on $I$, there exists a real constant $\delta$ such that
\Eq{W++}{ 
  W^{2,1}_{f,g}=\delta (W^{1,0}_{f,g})^3.
}
\end{enumerate}
}

An analogous proof of the following lemma has been introduced in \cite{PalZak19a}.
\Lem{DE}{Let $(f,g)\in\C_3(I)$. Then $f',g'$ form a fundamental system of solutions of the second-order homogeneous linear differential equation
\Eq{*}{
  W^{2,1}_{f,g}Y''=W^{3,1}_{f,g}Y'-W^{3,2}_{f,g}Y.
}}

Given an at most second-degree polynomial $P(u):=\alpha +\beta u+\gamma u^2$, where $\alpha,\beta,\gamma\in\R$, we call the value $D_P:=\beta^2-4\alpha\gamma$ the \emph{discriminant of $P$}. 

\Lem{DP}{If $P$ is an at most second-degree polynomial, then $D_P=(P')^2-2P''P$.}

\section{Main results}

For the proof of our main result we will need the following lemma, which describes an important property of pairs of functions belonging to the regularity class $\C_1(I)$. 

\Lem{inj}{If $(f,g)\in\C_1(I)$, then the mapping $(f,g):I\to\R^2$ is injective.}

\begin{proof}
To the contrary, assume that there exist $x<y$ in $I$ such that 
\Eq{*}{
  (f,g)(x)=(f,g)(y)=(p,q).
}
If, for all $z\in\,]x,y[\,$, the equality $(f,g)(z)=(p,q)$ holds, then $g'=0$ on $\,]x,y[\,$, contradicting that $g'$ is nonvanishing on $I$. Thus, there exists an element $z\in\,]x,y[\,$ such that $(f,g)(z)\neq(p,q)$. Applying the Cauchy Mean Value Theorem on the intervals $[x,z]$ and $[z,y]$, we can find two elements $u$ and $v$ with $x<u<z<v<y$ such that
\Eq{*}{
  f'(u)(g(z)-g(x))=g'(u)(f(z)-f(x)) \qquad\mbox{and}\qquad
  f'(v)(g(z)-g(y))=g'(v)(f(z)-f(y)).
}
Therefore, the vector $(\alpha,\beta):=(f(z)-p,g(z)-q)$ is a nontrivial solution of the following system of linear equations
\Eq{*}{
 f'(u)\beta-g'(u)\alpha=0 \qquad\mbox{and}\qquad
  f'(v)\beta-g'(v)\alpha=0.
}
Consequently, the determinant of this system must be zero, i.e., $f'(u)g'(v)=f'(v)g'(u)$. Dividing this equation by $g'(u)g'(v)$ side by side, it follows that $(f'/g')(u)=(f'/g')(v)$. On the other hand, our assumption implies that $f'/g'$ is strictly monotone, hence $u=v$, which contradicts $u<z<v$.
\end{proof}

\Lem{SM}{Let $(f,g)\in\C_1(I)$. Then $\CC_{f,g}$ is a symmetric, continuous and strictly monotone mean on $I$.}

\begin{proof} The symmetry and continuity are easy consequences of the definition of the Cauchy means. To prove the strict monotonicity in the first variable, let $x,y,z\in I$ with $x<y$. In the proof of the inequality
\Eq{xyz}{
  \CC_{f,g}(x,z)<\CC_{f,g}(y,z),
}
we assume that $f'/g'$ is strictly increasing, the other possibility is completely similar. If $z\in[x,y]$, then \eq{xyz} is a consequence of the strict mean property of $\CC_{f,g}$ because then
\Eq{*}{
  \CC_{f,g}(x,z)\leq z\leq \CC_{f,g}(y,z)
}
and one of the inequalities must be strict. Thus, we also may assume that $z\notin[x,y]$, that is, either $z<x$ or $y<z$. In these subcases \eq{xyz} is equivalent to
\Eq{*}{
  \frac{f(x)-f(z)}{g(x)-g(z)}<\frac{f(y)-f(z)}{g(y)-g(z)}.
}
Using that $g'$ is not vanishing, we have that $g$ is strictly monotone, therefore, the product of the denominators is positive in both subcases. Hence, the above inequality can be rewritten as
\Eq{*}{
  (f(x)-f(z))(g(y)-g(z))<(f(y)-f(z))(g(x)-g(z)).
}
This inequality is equivalent to
\Eq{*}{
  (f(x)-f(z))(g(y)-g(x))<(f(y)-f(x))(g(x)-g(z)).
}
Observe that in the first subcase $z<x$ the strict monotonicity of $g$ implies 
\Eq{*}{
  \frac{f(x)-f(z)}{g(x)-g(z)}<\frac{f(y)-f(x)}{g(y)-g(x)}.
}
Now, applying that $f'/g'$ is strictly increasing, the above inequality transforms to
\Eq{*}{
  \CC_{f,g}(x,z)<\CC_{f,g}(y,x).
}
This last inequality is seen to be true because, by the strict mean property of Cauchy means, $x$ separates the two sides. Hence \eq{xyz} holds as well. In the second subcase $y<z$, the proof is analogous.
\end{proof}

\Thm{M1}{
Let $(f,g)\in\C_1(I)$. Then the following statements are equivalent.
\begin{enumerate}[(i)]
\item There exists a continuous strictly monotone function $h:I\to\R$ such that
\Eq{CA}{
  \CC_{f,g}(x,y)=\A_h(x,y)\qquad(x,y\in I).
}
\item The mean $\CC_{f,g}$ is bisymmetric, i.e., it satisfies the following functional equation
\Eq{*}{
  \CC_{f,g}(\CC_{f,g}(x,y),\CC_{f,g}(u,v))
  =\CC_{f,g}(\CC_{f,g}(x,u),\CC_{f,g}(y,v)) \qquad(x,y,u,v\in I).
}
\item There exist real constants $\alpha,\beta,\gamma,\delta,\varepsilon,\eta$ with $(\alpha,\beta,\gamma)\neq(0,0,0)$ such that 
\Eq{CS}{
\alpha f^2+\beta fg+\gamma g^2+\delta f+\varepsilon g+\eta=0.
}
\item Provided that $(f,g)\in\C_2(I)$, there exist real constants $a,b, c$ such that
\Eq{QE}{
af'^2+bf'g'+cg'^2=(W^{2,1}_{f,g})^\frac{2}{3}.
}
\item Provided that $(f,g)\in\C_2(I)$,
\Eq{BA}{
  \B_{F,G}(x,y)=\A_h(x,y)\qquad(x,y\in I),
}
where $F:=f'/h'$, $G:=g'/h'$ and $h:=\int (W^{2,1}_{f,g})^\frac{1}{3}$.
\item Provided that $(f,g)\in\C_2(I)$, equation \eq{CA} holds with $h:=\int (W^{2,1}_{f,g})^\frac{1}{3}$.
\item Provided that $(f,g)\in\C_4(I)$, the expression
\Eq{exp}{
 \frac{3W^{4,1}_{f,g}+12W^{3,2}_{f,g}}{\Big(W^{2,1}_{f,g}\Big)^{\frac53}}
 -5\frac{\Big(W^{3,1}_{f,g}\Big)^2}{\Big(W^{2,1}_{f,g}\Big)^{\frac83}}\quad\mbox{is constant.}
}
\end{enumerate}
}

\begin{proof}
If $\CC_{f,g}$ is a quasiarithmetic mean, i.e., \eq{CA} holds with some strictly monotone and continuous function $h:I\to\R$, then, for all $x,y,u,v\in I$,
\Eq{*}{
  \CC_{f,g}&(\CC_{f,g}(x,y),\CC_{f,g}(u,v))
  =\A_h(\A_h(x,y),\A_h(u,v))
  =h^{-1}\bigg(\frac{h(\A_h(x,y))+h(\A_h(u,v))}{2}\bigg)\\
  &=h^{-1}\bigg(\frac{h(x)+h(y)+h(u)+h(v)}{4}\bigg) 
  =h^{-1}\bigg(\frac{h(x)+h(u)+h(y)+h(v)}{4}\bigg)\\
  &=h^{-1}\bigg(\frac{h(\A_h(x,u))+h(\A_h(y,v))}{2}\bigg)
  =\A_h(\A_h(x,u),\A_h(y,v)) 
  =\CC_{f,g}(\CC_{f,g}(x,u),\CC_{f,g}(y,v)).
}
The implication (ii)$\Rightarrow$(i) follows from Aczél's celebrated theorem \cite{Acz48a} (cf. \cite{Acz66}, \cite{AczDho89}) which says that every two-variable symmetric, continuous and strictly monotone mean which fulfils the bisymmetry property has to be a two-variable (symmetric) quasiarithmetic mean.
In our case, by \lem{SM}, $\CC_{f,g}$ is a symmetric, continuous and strictly monotone mean. Therefore, the result of Aczél directly applies.

In the next step prove that the assertions (i) and (iii) are equivalent. Assume that assertion (i) holds, i.e., there exists a continuous strictly monotone function $h:I\to\R$ such that \eq{CA} is valid. Applying \thm{0}, it follows that $h$ is differentiable with nonvanishing first derivative such that the equivalence \eq{equ} holds. Consequently, there exist real constants $a,b,c,d$ with $ad\neq bc$ such that
\Eq{SC-}{
  h'\cdot S_p\circ h=af'+bg' \qquad\mbox{and}\qquad h'\cdot C_p\circ h=cf'+dg'.
}
We consider two cases when we integrate these identities side by side.
If $p\neq0$, then we have the formulas
\Eq{*}{
  \int S_p=\frac{\sgn(p)}{\sqrt{|p|}}C_p,\qquad \int C_p=\frac{1}{\sqrt{|p|}}S_p.
}
Therefore from \eq{SC-}, it follows that there exist real constants $\lambda,\mu$ such that
\Eq{SC}{
C_p\circ h=a_1f+b_1g+\lambda \qquad\mbox{and}\qquad S_p\circ h=c_1f+d_1g+\mu
}
where $a_1:=\frac{\sqrt{|p|}}{\sgn(p)}a$, $b_1:=\frac{\sqrt{|p|}}{\sgn(p)}b$,
$c_1:=\sqrt{|p|}c$, and $d_1:=\sqrt{|p|}d$.
Using the well-known identities of trigonometric and hyperbolic functions, we have
\Eq{*}{
C^2_p-\sgn(p)S^2_p=1
}
holds on $\R$ and $C^2_p\circ h-\sgn(p)S^2_p\circ h=1$ is valid on $I$. Consequently, we obtain
\Eq{*}{
(a_1f+b_1g+\lambda)^2-\sgn(p)(c_1f+d_1g+\mu)^2=1,
}
and hence equation \eq{CS}  holds with the following constants
\Eq{*}{
\alpha&:=a_1^2-\sgn(p)c_1^2,\quad & \beta&:=2a_1b_1-2\sgn(p)c_1d_1,\quad & \gamma&:=b_1^2-\sgn(p)d_1^2,\\
\delta&:=2\lambda a_1-2\sgn(p)\mu c_1,\quad &
\varepsilon&:=2\lambda b_1-2\sgn(p)\mu d_1, \quad &
\eta&:=\lambda^2-\sgn(p)\mu^2-1.
}
To the contrary assume that $(\alpha,\beta,\gamma)=(0,0,0)$. If $p<0$, then these equalities imply that $a_1=b_1=c_1=d_1=0$, which yields $a=b=c=d=0$ contradicting $ad\neq bc$. In the case $p>0$, $(\alpha,\beta,\gamma)=(0,0,0)$ implies that $a_1^2=c_1^2$, $a_1b_1=c_1d_1$, and $b_1^2=d_1^2$. If $c_1=0$, then $a_1=0$ and hence
$ad=0=bc$, a contradiction. If $c_1\neq0$, then
\Eq{*}{
  a_1d_1=a_1\frac{a_1b_1}{c_1}=\frac{a_1^2b_1}{c_1}=\frac{c_1^2b_1}{c_1}=b_1c_1,
}
which again contradicts $ad\neq bc$. 

In the case $p=0$, the integration of the equalities \eq{SC-} yields the existence of constants $\lambda,\mu$ such that
\Eq{*}{
  \frac12 h^2=af+bg+\lambda \qquad\mbox{and}\qquad
  h=cf+dg+\mu.
}
Therefore,
\Eq{*}{
  \frac12 (cf+dg+\mu)^2=af+bg+\lambda.
}
Thus assertion (iii) is valid with the following constants
\Eq{*}{
\alpha&:=c^2,\quad & \beta&:=2cd,\quad & \gamma&:=d^2,\\
\delta&:=2\mu c-2a,\quad &
\varepsilon&:=2\mu d-2b, \quad &
\eta&:=\mu^2-2\lambda.
}
On the contrary suppose that $(\alpha,\beta,\gamma)=(0,0,0)$ which leads to $c=d=0$ contradicting $ad\neq bc$. Thus, we have shown that $(\alpha,\beta,\gamma)\neq(0,0,0)$ holds in all cases.

Now we prove that assertion (iii) implies (i). Consider the quadratic curve
\Eq{*}{
q:=\big\{(x,y)\in\R^2\mid\alpha x^2+\beta xy+\gamma y^2+\delta x+\varepsilon y+\eta=0\big\}.
}
By assumption (iii), this curve covers the range of the map $(f,g):I\to\R^2$. Therefore, $q$ cannot be empty or a single point. We are going to show that, in fact, this curve can only be either an ellipse, or a hyperbola, or a parabola.

There are three remaining degenerate cases concerning the curve $q$: 
\begin{enumerate}[(A)]
 \item $q$ is a straight line;
 \item $q$ is the union of two parallel lines;
 \item $q$ is the union of two intersecting lines.
\end{enumerate}
We prove by contradiction that none of these cases can happen. Assuming (A), (B), or (C), first we show that the range of $(f,g)$ is covered by one straight line. This is obvious in the case (A). In the case (B), the continuity of $(f,g)$ implies that its range is connected, hence it has to be contained in one of the parallel lines. Finally assume case (C), which implies that the curve $(f,g)$ is covered by two intersecting lines whose tangent unit vectors are denoted by $u$ and $v$. Since $g'$ is nowhere vanishing, thus the tangent vector field $(f',g')$ is also nowhere vanishing. On the other hand, this vector field is everywhere parallel either to $u$ or to $v$. Hence, the continuity of $(f',g')$ implies that it is everywhere parallel to one of them. This implies that the curve $(f,g)$ is covered by one of the lines.

Thus, we have proved that there exist three constants $\delta,\varepsilon,\eta\in\R$ with $(\delta,\varepsilon)\neq(0,0)$ such that
\Eq{*}{
  \delta f+\varepsilon g+\eta=0
}
holds on $I$. Differentiating this equality and dividing by $g'$, it follows that
\Eq{*}{
  \delta\frac{f'}{g'}=-\varepsilon,
}
which contradicts the strict monotonicity of $f'/g'$. This final contradiction yields that none of the cases (A), (B), or (C) can happen, and hence, $q$ can only be a nondegenerate conic section.

By elementary linear algebra, it follows that there exist six constants $a,b,c,d,\lambda,\mu\in\R$ with $ad\neq bc$ and two functions $\psi,\chi:I\to\R$ such that
\Eq{fg}{
\Big(\begin{array}{c}f\\g\end{array}\Big)=
\Big(\begin{array}{cc}a&b\\c&d\end{array}\Big)\Big(\begin{array}{c}\psi\\\chi\end{array}\Big)
+\Big(\begin{array}{c}\mu\\\lambda\end{array}\Big),
}
where $\psi,\chi$ satisfy one of the following equations:
\Eq{ehp}{
\psi^2+\chi^2=1, \qquad 
\psi^2-\chi^2=1, \qquad\mbox{and}\qquad
\psi=\chi^2.
}
Differentiating \eq{fg}, we obtain
\Eq{f'g'}{
\Big(\begin{array}{c}f'\\g'\end{array}\Big)=
\Big(\begin{array}{cc}a&b\\c&d\end{array}\Big)\Big(\begin{array}{c}\psi'\\\chi'\end{array}\Big).
}
According to \thm{0}, in all three cases we have to show that there is a number $p\in\R$ and a differentiable function $h\colon I\to\R$ with nonvanishing first derivative such that \eq{equ} holds.

First suppose that $(\psi,\chi)$ satisfies the first equation in \eq{ehp}. As we have seen in \lem{inj}, the map $(f,g)$ is injective. Therefore, the equality \eq{fg} implies that $(\psi,\chi)$ is also an injective map whose range is a subset of the unit circle by the first equation in \eq{ehp}. By the continuity of this map, we get that the range $R$ of $(\psi,\chi)$ is an open connected proper subset of the unit circle $S$. Then there exist $-2\pi<u<v<2\pi$ such that the range of the map $(\cos,\sin)$ restricted to the interval $\,]u,v[\,$ equals $R$. Define the transformation $T:\R_+\times \,]u,v[\,\to\R^2$ by 
\Eq{*}{
  T(r,t):=(r\cos t,r\sin t).
}
Then $T$ is an injective differentiable map whose derivative matrix is nonsingular at every point in $\R_+\times \,]u,v[\,$, therefore, the inverse of $T$ is differentiable by the inverse function theorem. Finally, define $h:I\to\,]u,v[\,$ as the second coordinate function of $T^{-1}\circ(\psi,\chi)$. Then $h$ is differentiable and the equalities $\psi=\cos \circ h$ and $\chi=\sin\circ h$ hold. We can calculate the derivative of $f$ and $g$:
\Eq{*}{
f'&=(-a\sin\circ h+b\cos\circ h)h'=(-a S_{-1}\circ h+b C_{-1}\circ h)h', \\
g'&=(-c\sin\circ h+d\cos\circ h)h'=(-c S_{-1}\circ h+d C_{-1}\circ h)h'.
}
Since we assumed that $g'$ never vanishes, from this it also follows that $h'$ never vanishes. Thus the condition \eq{equ} of \thm{0} is satisfied with $p=-1$, and hence (i) holds.

Secondly, assume that $(\psi,\chi)$ fulfils the second equation in \eq{ehp}. Then define \mbox{$h\colon I\to\R$} by $h:=\sinh^{-1}\circ\chi$. Therefore $h$ is differentiable, $\chi=\sinh\circ h$, and, by the second equation in \eq{ehp}, we have that $\psi=\pm\cosh\circ h$. Thus, using \eq{f'g'}, for the derivatives of $f$ and $g$, we obtain
\Eq{*}{
f'&=(\pm a\sinh\circ h+b\cosh\circ h)h'=(\pm a S_1\circ h+b C_1\circ h)h', \\
g'&=(\pm c\sinh\circ h+d\cosh\circ h)h'=(\pm c S_1\circ h+d C_1\circ h)h'.
}
Again we can see that $h'$ is nonvanishing. Thus the condition \eq{equ} of \thm{0} is now satisfied with $p=1$, and consequently (i) is valid.

Finally, if the third equality of \eq{ehp} holds, then let $h:=\chi$, which is now automatically differentiable. Then $\psi=\chi^2=h^2$, and now using \eq{f'g'} we can calculate the derivatives of $f$ and $g$:
\Eq{*}{
f'&=(2a h+b)h'=(2aS_0\circ h+bC_0\circ h)h',\\
g'&=(2c h+d)h'=(2cS_0\circ h+dC_0\circ h)h'.
}
The second equation implies that $h'$ is nonvanishing. Therefore the relation \eq{equ} of \thm{0} is again satisfied by $p=0$, which yields condition (i).

To prove the implication (iii)$\Rightarrow$(iv), assume that $(f,g)\in\C_2(I)$. Assume that (iii) holds for some $\alpha,\beta,\gamma,\delta,\varepsilon,\eta\in\R$ with $(\alpha,\beta,\gamma)\neq(0,0,0)$. Differentiating \eq{CS}, we get 
\Eq{*}{
2\alpha f'f+\beta(f'g+fg')+2\gamma g'g+\delta f'+\varepsilon g'=0.
}
Denote $\varphi:=f'/g'$. Then, from the assumption $(f,g)\in\C_2(I)$ it follows that $\varphi$ is continuously differentiable and $\varphi'$ is nowhere zero on $I$. Now replacing $f'$ by $\varphi g'$, we obtain
\Eq{CS'}{
\varphi (2\alpha f+\beta g+\delta)+\beta f+2\gamma g+\varepsilon=0.
}
Differentiating this equation and then replacing $f'$ by $\varphi g'$, we arrive at
\Eq{CS''}{
\varphi'(2\alpha f+\beta g+\delta)+2g'(\alpha\varphi^2+\beta\varphi+\gamma)=0.
}
This implies that
\Eq{CS'+}{
  2\alpha f+\beta g+\delta=-\frac{2g'}{\varphi'}(\alpha\varphi^2+\beta\varphi+\gamma).
}
Observe that the last factor of the right hand side is a nontrivial at most second degree polynomial of $\varphi$. The function $\varphi$ is strictly monotone, therefore, the right hand side and consequently the left hand side of \eq{CS'+} can have at most two distinct zeros whose set will be denoted by $Z$.

Then, by \eq{CS'}, on the set $I\setminus Z$, we can write 
\Eq{*}{
\varphi=-\frac{\beta f+2\gamma g+\varepsilon}{2\alpha f+\beta g+\delta}.
}
This equality, combined with $(f,g)\in\C_2(I)$, implies that $\varphi$ is twice differentiable on $I\setminus Z$. 

Differentiating \eq{CS''} and replacing $f'$ by $\varphi g'$, on the set $I\setminus Z$, we get
\Eq{CS'''}{
\varphi''(2\alpha f+\beta g+\delta)+2g''(\alpha\varphi^2+\beta\varphi+\gamma)+3g'(2\alpha\varphi'\varphi+\beta\varphi')=0.
}
Using \eq{CS'+}, the above equality reduces to
\Eq{CS'''+}{
-2\frac{\varphi''}{\varphi'}g'(\alpha\varphi^2+\beta\varphi+\gamma)+2g''(\alpha\varphi^2+\beta\varphi+\gamma)+3g'(2\alpha\varphi'\varphi+\beta\varphi')=0.
}
Then, dividing this equation side by side by $3g'(\alpha\varphi^2+\beta\varphi+\gamma)$ on the set $I\setminus Z$, we obtain 
\Eq{*}{
\frac{2}{3}\frac{g''}{g'}-\frac{2}{3}\frac{\varphi''}{\varphi'}+\frac{2\alpha\varphi'\varphi+\beta\varphi'}{\alpha\varphi^2+\beta\varphi+\gamma}=0.
}
Integrating both sides, it follows that
\Eq{*}{
\frac23\ln|g'|-\frac23\ln|\varphi'|+\ln|\alpha\varphi^2+\beta\varphi+\gamma|
}
equals a constant on each component of $I\setminus Z$. Therefore,
\Eq{*}{
\Phi:=\Big(\frac{g'}{\varphi'}\Big)^\frac{2}{3}(\alpha\varphi^2+\beta\varphi+\gamma)
}
equals a nonzero constant on each component of $I\setminus Z$. On the other hand, $\Phi$ is continuous on $I$, the set $Z$ contains at most two points, consequently $\Phi$ is identically equal to a nonzero constant $\zeta$ on $I$.
Combining this result with equalities $\varphi=f'/g'$ and
\Eq{W21}{
\varphi'=\bigg(\frac{f'}{g'}\bigg)'=\frac{f''g'-g''f'}{g'^2}
  =\frac{W_{f,g}^{2,1}}{g'^2},
}  
we get assertion (iv) with constants $a:=\alpha/\zeta$, $b:=\beta/\zeta$ and $c:=\gamma/\zeta$.
 
To prove the implication (iv)$\Rightarrow$(v), assume that $(f,g)\in\C_2(I)$. If (iv) holds, then there exist real constants $a,b,c$ such that \eq{QE} is valid. Denote $F:=f'/h'$, $G:=g'/h'$ and $h:=\int (W^{2,1}_{f,g})^\frac{1}{3}$, then equation \eq{QE} reduces to
\Eq{QE+}{
aF^2+bFG+cG^2=1,
}
where $G$ is nowhere zero on $I$ and $F/G=f'/g'$ is strictly monotone and continuous on $I$. Applying implication (iii)$\Rightarrow$(i) of \thm{0+}, we conclude that assertion (v) holds.

Assume now that assertion (v) is valid, i.e., $(f,g)\in\C_2(I)$ and the functional equation \eq{BA} satisfied with $F:=f'/h'$, $G:=g'/h'$ and $h:=\int (W^{2,1}_{f,g})^\frac{1}{3}$. Applying implication (i)$\Rightarrow$(ii) of \thm{0+}, we get
\Eq{*}{
(F,G)\sim(S_p\circ h,C_p\circ h),
}
or equivalently, 
\Eq{*}{
(f',g')\sim(h'\cdot S_p\circ h,h'\cdot C_p\circ h). 
}
Therefore, using \thm{0}, we get assertion (vi). The implication (vi)$\Rightarrow$(i) is obvious. Hence all the assertions from (i) to (vi) are equivalent provided that $(f,g)\in\C_2(I)$.

To prove the implication (iv)$\Rightarrow$(vii), assume that $(f,g)\in\C_4(I)$. If (iv) holds, then there exist real constants $a,b,c$ such that equation \eq{QE} is valid. Denoting $\varphi:=f'/g'$ and replacing $f'$ by $\varphi g'$ in \eq{QE}, we obtain
\Eq{P}{
P\circ\varphi=\frac{(W^{2,1}_{f,g})^\frac{2}{3}}{g'^2}, 
}
where $P$ is an at most second-degree polynomial. Differentiating equation \eq{P}, we arrive at
\Eq{*}{
(P'\circ\varphi)\varphi'=\frac{2}{3g'^2}(W^{2,1}_{f,g})^\frac{-1}{3}W^{3,1}_{f,g}
-2\frac{g''}{g'^3}(W^{2,1}_{f,g})^\frac{2}{3}.
}
Using identity \eq{W21}, this equation reduces to
\Eq{P'}{
P'\circ\varphi=\frac{2}{3}(W^{2,1}_{f,g})^\frac{-4}{3}W^{3,1}_{f,g}
-2\frac{g''}{g'}(W^{2,1}_{f,g})^\frac{-1}{3}.
}
Differentiating equation \eq{P'}, we get
\Eq{*}{
(P''\circ\varphi)\varphi'=&-\frac{8}{9}(W^{2,1}_{f,g})^\frac{-7}{3}(W^{3,1}_{f,g})^2
+\frac{2}{3}(W^{2,1}_{f,g})^\frac{-4}{3}(W^{3,2}_{f,g}+W^{4,1}_{f,g})\\
&-2\frac{g'''g'-g''^2}{g'^2}(W^{2,1}_{f,g})^\frac{-1}{3}
+\frac23\frac{g''}{g'}(W^{2,1}_{f,g})^\frac{-4}{3}W^{3,1}_{f,g}.
}
Again using identity \eq{W21}, this equation simplifies to 
\Eq{P''}{
P''\circ\varphi=&-\frac{8}{9\varphi'}(W^{2,1}_{f,g})^\frac{-7}{3}(W^{3,1}_{f,g})^2
+\frac{2}{3\varphi'}(W^{2,1}_{f,g})^\frac{-4}{3}(W^{3,2}_{f,g}+W^{4,1}_{f,g})\\
&-2(g'''g'-g''^2)(W^{2,1}_{f,g})^\frac{-4}{3}
+\frac23\frac{g''}{g'\varphi'}(W^{2,1}_{f,g})^\frac{-4}{3}W^{3,1}_{f,g}.
}
Therefore, using \lem{DP} and identity \eq{W21}, we obtain
\Eq{*}{
D_P=&(P'\circ\varphi)^2-2(P''\circ\varphi)(P\circ\varphi)\\
=&\frac{4}{9}(W^{2,1}_{f,g})^\frac{-8}{3}(W^{3,1}_{f,g})^2
-\frac{8}{3}\frac{g''}{g'}(W^{2,1}_{f,g})^\frac{-5}{3}W^{3,1}_{f,g}
+4\frac{g''^2}{g'^2}(W^{2,1}_{f,g})^\frac{-2}{3}\\
&+\frac{16}{9}(W^{2,1}_{f,g})^\frac{-8}{3}(W^{3,1}_{f,g})^2
-\frac{4}{3}(W^{2,1}_{f,g})^\frac{-5}{3}(W^{3,2}_{f,g}+W^{4,1}_{f,g})\\
&+4\frac{g'''g'-g''^2}{g'^2}(W^{2,1}_{f,g})^\frac{-2}{3}
-\frac43\frac{g''}{g'}(W^{2,1}_{f,g})^\frac{-5}{3}W^{3,1}_{f,g},
}
or equivalently,
\Eq{*}{
D_P=\frac{20}{9}(W^{2,1}_{f,g})^\frac{-8}{3}(W^{3,1}_{f,g})^2
-\frac{4}{3}(W^{2,1}_{f,g})^\frac{-5}{3}(W^{3,2}_{f,g}+W^{4,1}_{f,g})
-4(W^{2,1}_{f,g})^\frac{-5}{3}\frac{g''W^{3,1}_{f,g}-g'''W^{2,1}_{f,g}}{g'}.
}
It is easy to check that $\dfrac{g''W^{3,1}_{f,g}-g'''W^{2,1}_{f,g}}{g'}=W^{3,2}_{f,g}$. Therefore, we get the expression \eq{exp} is equal to $-\frac94 D_p$. Hence assertion (vii) holds. 

To complete the proof of the theorem it suffices to prove the implication (vii)$\Rightarrow$(i) in the class $\C_4(I)$. Assume that (vii) holds, i.e., the expression in \eq{exp} is equal to constant. Let $Y\in\{f',g'\}$, using \lem{DE}, it follows that $Y$ is a solution of the following second-order homogeneous linear differential equation
\Eq{DE}{
  W^{2,1}_{f,g}Y''=W^{3,1}_{f,g}Y'-W^{3,2}_{f,g}Y.
}
Now, denote $h:=\int (W^{2,1}_{f,g})^\frac{1}{3}$. It follows that $h$ is three times differentiable strictly monotone with a nonvanishing first derivative. Therefore, its inverse is also three times differentiable. Define the function $Z:h(I)\to\R$ by $Z:=\frac{1}{h'\circ h^{-1}}(Y\circ h^{-1})$. Consequently, $Z$ is a twice differentiable function and we have $Y=h'(Z\circ h)$. Differentiating $Y$ once and twice, we get
\Eq{*}{
Y'=h''(Z\circ h)+h'^2(Z'\circ h)\qquad\mbox{and}\qquad Y''=h'''(Z\circ h)+3h''h'(Z'\circ h)+h'^3(Z''\circ h).
}
However we have,
\Eq{*}{
h'=(W^{2,1}_{f,g})^\frac{1}{3},\quad h''=\frac13(W^{2,1}_{f,g})^\frac{-2}{3}W^{3,1}_{f,g},\quad 
h'''=-\frac{2}{9}(W^{2,1}_{f,g})^\frac{-5}{3}(W^{3,1}_{f,g})^2+\frac13(W^{2,1}_{f,g})^\frac{-2}{3}(W^{4,1}_{f,g}+W^{3,2}_{f,g}). 
}
Applying these identities and \eq{DE}, we arrive at
\Eq{*}{
\big(&-\frac{2}{9}(W^{2,1}_{f,g})^\frac{-2}{3}(W^{3,1}_{f,g})^2+\frac13(W^{2,1}_{f,g})^\frac{1}{3}(W^{4,1}_{f,g}+W^{3,2}_{f,g})\big)(Z\circ h)+(W^{2,1}_{f,g})^\frac{2}{3}W^{3,1}_{f,g}(Z'\circ h)+(W^{2,1}_{f,g})^2(Z''\circ h)\\
&=\frac13(W^{2,1}_{f,g})^\frac{-2}{3}(W^{3,1}_{f,g})^2(Z\circ h)+(W^{2,1}_{f,g})^\frac{2}{3}W^{3,1}_{f,g}(Z'\circ h)
-(W^{2,1}_{f,g})^\frac{1}{3}W^{3,2}_{f,g}(Z\circ h).
}
This simplifies to the identity
\Eq{*}{
Z''\circ h=-\frac19\Bigg(\frac{3W^{4,1}_{f,g}+12W^{3,2}_{f,g}}{\Big(W^{2,1}_{f,g}\Big)^{\frac53}}
 -5\frac{\Big(W^{3,1}_{f,g}\Big)^2}{\Big(W^{2,1}_{f,g}\Big)^{\frac83}}\Bigg) (Z\circ h).
}
Therefore, using assertion (vii), there exists real constant $p$ such that $Z''\circ h=pZ\circ h$ is valid on $I$. Thus, it follows that 
\Eq{Z''}{
Z''=pZ 
}
holds on $h(I)$. Thus, $Z:=\frac{1}{h'\circ h^{-1}}(f'\circ h^{-1})$ and $Z:=\frac{1}{h'\circ h^{-1}}(g'\circ h^{-1})$ are solutions to this second-order homogeneous linear differential equation. On the other hand $(S_p,C_p)$ forms a fundamental solution system for \eq{Z''}. Consequently,
\Eq{*}{
\bigg(\frac{1}{h'\circ h^{-1}}(f'\circ h^{-1}),\frac{1}{h'\circ h^{-1}}(g'\circ h^{-1})\bigg)\sim(S_{p},C_{p}).
}
Thus, the relation \eq{equ} is satisfied so we conclude that the assertion (i) holds.
\end{proof}

\def\MR#1{}


\begin{thebibliography}{10}

\bibitem{Acz48a}
J.~Aczél, \emph{{On mean values}}, Bull. Amer. Math. Soc. \textbf{54} (1948),
  392–400. \MR{9,501h}

\bibitem{Acz66}
J.~Aczél, \emph{{{L}ectures on {F}unctional {E}quations and {T}heir
  {A}pplications}}, {Mathematics in Science and Engineering}, vol.~19, Academic
  Press, New York–London, 1966. \MR{34 \#8020}

\bibitem{AczDho89}
J.~Aczél and J.~Dhombres, \emph{{Functional {E}quations in {S}everal
  {V}ariables}}, Cambridge University Press, Cambridge, 1989, With applications
  to mathematics, information theory and to the natural and social sciences.
  \MR{90h:39001}

\bibitem{Baj58}
M.~Bajraktarević, \emph{{Sur une équation fonctionnelle aux valeurs
  moyennes}}, Glasnik Mat.-Fiz. Astronom. Društvo Mat. Fiz. Hrvatske Ser. II
  \textbf{13} (1958), 243–248. \MR{23 \#A442}

\bibitem{DarMakPal04}
Z.~Daróczy, Gy. Maksa, and Zs. Páles, \emph{{On two-variable means with
  variable weights}}, Aequationes Math. \textbf{67} (2004), no.~1-2, 154–159.
  \MR{2005a:39043}

\bibitem{HarLitPol34}
G.~H. Hardy, J.~E. Littlewood, and G.~Pólya, \emph{{Inequalities}}, Cambridge
  University Press, Cambridge, 1934, (first edition), 1952 (second edition).
  \MR{13,727e}

\bibitem{KisPal18}
T.~Kiss and Zs. Páles, \emph{{On a functional equation related to two-variable
  weighted quasi-arithmetic means}}, J. Difference Equ. Appl. \textbf{24}
  (2018), no.~1, 107–126. \MR{3750533}

\bibitem{KisPal19}
T.~Kiss and Zs. Páles, \emph{{On a functional equation related to two-variable Cauchy
  means}}, Math. Inequal. Appl. (2019), to appear.

\bibitem{LeaSho84}
E.~Leach and M.~Sholander, \emph{{Multivariable extended mean values}}, J.
  Math. Anal. Appl. \textbf{104} (1984), 390–407.

\bibitem{Los00a}
L.~Losonczi, \emph{{Equality of {C}auchy mean values}}, Publ. Math. Debrecen
  \textbf{57} (2000), 217–230. \MR{2001g:39051}

\bibitem{PalZak19c}
Zs. Páles and A.~Zakaria, \emph{{On the equality of Bajraktarević means to
  quasi-arithmetic means}},  (2019), submitted.

\bibitem{PalZak19a}
Zs. Páles and A.~Zakaria, \emph{{On the invariance equation for two-variable weighted
  nonsymmetric Bajraktarević means}}, Aequationes Math. \textbf{93} (2019),
  no.~1-2, 37–57. \MR{3919421}

\end{thebibliography}

\end{document}